\documentclass[amsppt,11pt]{amsart}
\usepackage{curves}
\usepackage{amsmath}

\newtheorem {theorem}{Theorem}[section]

\newtheorem {proposition}[theorem]{Proposition}
\newtheorem {lemma}[theorem]{Lemma}

\newcounter{conjecture}
\newcounter{remark}

\newcommand{\eqnsection}{
   \renewcommand{\theequation}{\thesection.\arabic{equation}}
   \makeatletter
   \csname @addtoreset\endcsname{equation}{section}
   \makeatother}

\newcommand{\be}{{\begin{equation}}}
\newcommand{\ee}{{\end{equation}}}
\def \bt{\begin{theorem}}
\def \et{\end{theorem}}
\def \bea{\begin{eqnarray}}
\def \eea{\end{eqnarray}}
\def \bas{\begin{eqnarray*}}
\def \eas{\end{eqnarray*}}

\def \de{\delta}

\def \ep{\epsilon}
\def \epon{\ep_1}

\newcommand{\eps}{\varepsilon}
\newcommand{\bart}{{\bar{\theta}}}
\newcommand{\bars}{{\bar{\sigma}}}

\def \om{\omega}
\def \Om{\Omega}

\def \th{\theta}

\def \ff{\infty}
\def \wh{\widehat}

\def \rar{\rightarrow}

\newcommand{\ls}[1]
   {\dimen0=\fontdimen6\the\font \lineskip=#1\dimen0
\advance\lineskip.5\fontdimen5\the\font \advance\lineskip-\dimen0
\lineskiplimit=.9\lineskip \baselineskip=\lineskip
\advance\baselineskip\dimen0 \normallineskip\lineskip
\normallineskiplimit\lineskiplimit \normalbaselineskip\baselineskip
\ignorespaces }

\def \R{{\Bbb{R}}}

\def \Z{{\Bbb{Z}}}

\def \AA{{\mathcal A}}

\def \HH{{\mathcal H}}
\def \II{{\mathcal I}}

\def \TT{{\mathcal T}}

\def \VV{{\mathcal V}}

\def \({\left(}
\def \){\right)}

\def \nn{\nonumber}

\def \bc{\begin{center} }
\def \ec{\end{center} }
\def\Bbb{\mathbb}

\begin{document}

\eqnsection
\newcommand{\Ini}{{I_{n,i}}}
\newcommand{\reals}{{\Bbb{R}}}
\newcommand{\F}{{\mathcal F}}
\newcommand{\D}{{\mathcal D}}
\newcommand{\Fn}{{{\mathcal F}_n}}
\newcommand{\Gn}{{{\mathcal G}_n}}
\newcommand{\Hn}{{{\mathcal H}_n}}
\newcommand{\Fp}{{{\mathcal F}^p}}
\newcommand{\Gp}{{{\mathcal G}^p}}
\newcommand{\PPP}{{\mathbf P}}
\newcommand{\Pop}{{P\otimes \PPP}}
\newcommand{\hm}{\HH^\varphi}
\newcommand{\nuw}{{\nu^W}}
\newcommand{\ths}{{\theta^*}}
\newcommand{\beq}[1]{\begin{equation}\label{#1}}
\newcommand{\eeq}{\end{equation}}
\newcommand{\integers}{{\rm I\!N}}
\newcommand{\DT}{{D_{\Bbb T^2}}}
\newcommand{\pDT}{{\partial \DT}}
\newcommand{\E}{{\Bbb E}}
\newcommand{\te}{{\tilde{\delta}}}
\newcommand{\tI}{{\tilde{I}}}
\newcommand{\epn}{{\ep_n}}
\def\var{{\rm Var}}
\def\cov{{\rm Cov}}
\def\one{{\bf 1}}
\def\leb{{\mathcal L}eb}
\def\Ho{{\mbox{\sf H\"older}}}  
\def\thi{{\mbox{\sf Thick}}}
\def\cthi{{\mbox{\sf CThick}}}
\def\late{{\mbox{\sf Late}}}
\def\clate{{\mbox{\sf CLate}}}
\def\plate{{\mbox{\sf PLate}}}
\newcommand{\ffrac}[2]
  {\left( \frac{#1}{#2} \right)}
\newcommand{\calF}{{\mathcal F}}
\newcommand{\dfn}{\stackrel{\triangle}{=}}
\newcommand{\beqn}[1]{\begin{eqnarray}\label{#1}}
\newcommand{\eeqn}{\end{eqnarray}}
\newcommand{\oo}{\overline}
\newcommand{\uu}{\underline}
\newcommand{\bfcdot}{{\mbox{\boldmath$\cdot$}}}
\newcommand{\Var}{{\rm \,Var\,}}
\def\squarebox#1{\hbox to #1{\hfill\vbox to #1{\vfill}}}
\renewcommand{\qed}{\hspace*{\fill}
            \vbox{\hrule\hbox{\vrule\squarebox{.667em}\vrule}\hrule}\smallskip}
\newcommand{\half}{\frac{1}{2}\:}
\newcommand{\beaa}{\begin{eqnarray*}}
\newcommand{\eeaa}{\end{eqnarray*}}
\newcommand{\calK}{{\mathcal K}}
\def\dimm{{\overline{{\rm dim}}_{_{\rm M}}}}
\def\dimp{\dim_{_{\rm P}}}
\def\htaum{{\hat\tau}_m}
\def\htaumk{{\hat\tau}_{m,k}}
\def\htaumkj{{\hat\tau}_{m,k,j}}

\bibliographystyle{amsplain}

\title[Thick Points for the Cauchy Process]
{Thick Points for the Cauchy Process}

\author[Olivier Daviaud]
{Olivier Daviaud$^*$}

\date{June 9, 2003.
\newline\indent
$^*$Department of Mathematics, Stanford University, Stanford CA 94305, USA. Email: 
\indent odaviaud@stanford.edu. 
Research  partially supported by NSF grant \#DMS-0072331.
}

\begin{abstract}
\noindent
Let $\TT(x,\eps)$ denote the occupation measure of an interval of length $2\eps$ centered at $x$
by the Cauchy process run until it hits $(-\infty,-1]\cup [1,\infty)$. 
We prove that $\sup_{|x|\leq 1}\TT(x,\eps)/(\eps(\ln\eps)^2)\to 2/\pi$ a.s. as $\eps\to 0$. 
We also obtain the multifractal spectrum for thick points, i.e.
the Hausdorff dimension of the set of $\alpha$-thick points $x$ for which
$\lim_{\eps \to 0} \TT(x,\eps)/(\eps(\ln\eps)^2) = \alpha > 0$.
\end{abstract}
\maketitle

\section{Introduction}
Let $X=(X_t, t\geq 0)$ be a Cauchy process on the real line $\mathbb{R}$, that is a process starting at $0$, 
with stationary independent increments with the Cauchy  distribution:
$$P(X_{t+s}-X_t\in dx)=\frac{s\,dx}{\pi(s^2+x^2)},\quad s,t>0,\ x\in\mathbb{R}.$$ 
Next, let 
$$\mu_\bart^X(A):=\int_0^\bart 1_A(X_s)ds,$$
be the occupation measure of a measurable subset $A$ of $\reals$ by the Cauchy process run until $\bart:=\inf\{s:|X_s|\geq1\}$.
Let $I(x,\ep)$ denote the interval of radius $\ep$ centered at $x$.
Our first theorem follows:
\begin{theorem}
\label{theo-sup}
\begin{equation}
\label{eq-sup}
\lim_{\ep\to 0}\sup_{x\in \mathbb{R}}
\frac{\mu_\bart^X(I(x,\ep))}{\ep (\log \ep)^2}=2/\pi \hspace{.3in}a.s.
\end{equation}
\end{theorem}
Compare our result to the analogue of Ray's result (\cite{Ray}):
$$\limsup_{\ep\to 0}\frac{\mu_\bart^X(I(0,\ep))}{\ep\log\frac{1}{\ep}\log\log\log\frac{1}{\ep}}=c \hspace{.3in}a.s.$$
for some constant $0<c<\infty$. Indeed, this can be proved by slightly modifying the proof in \cite{Ray} ( 
using the system of excursions that we introduce in Section 4).

Next, it follows from the previous theorem (or more simply from \cite[Lemma 2.3]{PT}) that for almost all paths,
$$\limsup_{\eps\to 0}\frac{\log\mu_\bart^X(I(x,\ep))}{\log\ep}\leq 1$$
for \emph{all} points $x$ in the range $\{X_t\ |\,0\leq t\leq\bart\}$. On the other hand, this fact together with 
\cite[Chap. VIII Theorem 5]{Bertoin} and Fubini's theorem imply that
for $\PPP\times\leb$-almost every $(\om,t)$ in $\Om\times\{0\leq t\leq\bart\}$
$$\lim_{\eps\to 0}\frac{\log\mu_\bart^X(I(x,\ep))}{\log\ep}=1.$$
Hence, standard multifractal analysis must
be refined in order to obtain a non-degenerate dimension spectrum for thick points. This leads us to 
\begin{theorem}
\label{theo-thick}
For any $a\leq 2/\pi$,
\begin{equation}\label{eq-thick}
\dim\Big\{x\in \reals:\;
\lim_{\eps\to0} \frac{\mu_\bart^X (I(x,\eps))}{\eps
\left(\log \eps\right)^2}= a\Big\}= 1-a\pi/2
\hspace{.3in}a.s.
\end{equation}
\end{theorem}

The results obtained here are the analogues of those in \cite{DPRZ4}, when replacing the planar Brownian motion (which is a stable process 
of index $2$ in dimension $2$) by the Cauchy process (which is a stable process of index $1$ in dimension 1). Theorems \ref{theo-sup} and
\ref{theo-thick} answer the first part of open problem $(6)$ of that paper (also implicitly present in \cite{PT}), 
the second part being solved in \cite{DPRZ5}. Our work relies heavily on
the techniques developed in \cite{DPRZ4} and \cite{DPRZ5}, and therefore owes a substantial debt to these papers. 

The Cauchy process is a symmetric stable process of index $\alpha=1$. Results similar to those of this paper have been obtained for $\alpha<1$ 
(i.e. for transient symmetric stable processes in one dimension) in \cite{DPRZ3}. The case $\alpha>1$ is easier to study, since for such processes
there exists a bi-continuous local time (e.g. see \cite{Boylan}). Thus in that case Theorem \ref{theo-sup} would hold, but with a different
scaling (simply $\ep$) and $2/\pi$ would be replaced by a random variable (more precisely: the supremum of the local time). Thus our findings apply
only at the border of transience and recurrence of stable processes.

The main difficulty in obtaining results similar to those in \cite{DPRZ4} is that the Cauchy process is not continuous. Indeed, the proof of the 
lower bounds
in \cite{DPRZ4} relies on the idea that unusually high occupation measures in the neighborhood of a point $x$ are the result of an unusually high number 
of excursions of all scales around this point. But defining the notion of excursion is not clear when it comes to a non-continuous process.
Our proof avoids this problem essentially by working with the Brownian representation of the Cauchy process: up to a time-change, the Cauchy process
can be seen as the intersection of a two-dimensional Brownian motion and, say, the $x$-axis. Using this framework we obtain lower bounds by adapting the
strategy in \cite{DPRZ5}. The same strategy could be used to derive upper bound results. However, because of its independent interest we use the following proposition
as the key to our proof of the upper bounds:
\begin{proposition}\label{prop-GreenForCauchy}
Fix $r>0$, and let $\bart=\inf\{t: |X_t|\geq r\}$. For any bounded Borel measurable function 
$f:[-r,r]\to\reals$,
\begin{equation}\label{eq-GreenForCauchy}
\E^{x_0}\int_0^\bart f(X_s)ds=\int_{-r}^r f(x)G(x_0,x)dx
\end{equation}
where $G$ is given by 
\begin{equation}\label{eq:G}
G(x_0,x)=-(1/2\pi)\log\left|\frac{h(x/r)-h(x_0/r)}{1-h(x/r)h(x_0/r)}\right|,
\end{equation}
and
\begin{equation}\label{eq:H}
h(x)=\frac{\sqrt{\frac{1+x}{1-x}}-1}{\sqrt{\frac{1+x}{1-x}}+1}.
\end{equation}
\end{proposition}

\noindent
\textbf{Remarks.}

\noindent
$\bullet \;$ By the scaling property of the Cauchy process, for any deterministic $0<r<\infty$, Theorems \ref{theo-sup} and \ref{theo-thick}
still hold if we replace $\bart$ by $\bart_r=\inf\{t:\:|X_t|\geq r\}$. As a consequence, these results also hold if one replaces $\bart$ by any
deterministic $T<\infty$, or any almost surely finite stopping time.

\noindent
$\bullet \;$ In the course of our study, we will prove (see equation (\ref{eq-upperthick})) that almost surely
$$\dim\{x\in I(0,1): \lim_{\ep\to 0}
\frac{\mu_\bart^X(I(x,\ep))}{\ep (\log \ep)^2}\geq a\}\leq 1-a\pi/2.$$
Using this fact, Theorem \ref{theo-thick} still holds if in equation (\ref{eq-thick}) one replaces $\lim$ by $\limsup$ or $\liminf$,
and/or '$=a$' by '$\geq a$'.

\noindent
$\bullet \;$ Exactly as in \cite{DPRZ4}, one can obtain the following result for the coarse multi-fractal spectrum: for every $a<2/\pi$,
$$ \lim_{\eps \to 0} \frac{\log \leb (x: \mu_\bart^X(I(x,\eps))\geq
a
\eps
(\log \eps)^2)}{\log \eps} = a\frac{\pi}{2}\,, \quad a.s.$$

It is quite natural to consider also the discrete analogues of the results presented here. For example, let 
$(X_i)$ be a sequence of i.i.d. variables with distribution:
$$\PPP(X_i=n)=\frac{C}{1+n^2}\,,\,\,n\in\mathbb{Z}$$
where $C$ is a normalizing constant. Let $S_n=\sum_{1}^n X_i$,
$$
L_n^{S}(x):=\#\{i: S_i=x, \,0 \leq i\leq n\}\,,\,\,\,
$$
be the number of visits to $x\in \Z$ during the first $n$ steps of the walk
and 
$$
T_n^{X}:= \max_{x\in \Z^2} L_n^{X}(x)\,,
$$
its maximal value. Then we conjecture that there exists a constant $\alpha$ such that
$$
\lim_{n\to\infty} \frac{T_n^{X}}{(\log n)^2}=
\alpha\,,
\,\,\,a.s.$$
The source of the difficulty here is the absence of strong approximation theorems (such results were used to prove the Erd\H{o}s-Taylor conjecture
in \cite{DPRZ4}). Another integer-valued random variable for which we expect similar asymptotic results is the following one:
 if $X_n=(X_n^1,X_n^2)$ is a simple random walk in $\Z^2$, then we define 
$Y_n:=X^1_{t_n}$ where $t_n$ is the time at which $0$ is visited for the $n^{\textrm{th}}$ time by $X^2$. This is the discrete time analogue of the Brownian 
representation of the Cauchy process, so our techniques should apply here. More generally, we suspect the existence of similar results for
random variables in the domain of attraction of the Cauchy distribution.

In the next section, we prove Proposition \ref{prop-GreenForCauchy} using the Brownian representation of the Cauchy process
and the solution to some Dirichlet problem. In Section 3, we use this result to prove upper bounds for both theorems. 
In Section 4 we prove the lower bounds, using a well defined system of
excursions analogous to the one which appears in \cite{DPRZ5}. Finally Section 5 establishes the connection between occupation measure and excursions.

\section{Green function for the Cauchy process}
This section will be devoted to proving Proposition \ref{prop-GreenForCauchy}.
By a density argument, it is enough to prove the proposition for $f$ continuous with compact support in $(-r,r)$. 
The proof is based on the 
Brownian representation of the Cauchy process: if $(B^1,B^2)$ is a planar Brownian motion, and $\tau(t)$ is the inverse local time 
of $B^2$ at $0$, then $B^1_{\tau(t)}$ is a Cauchy process. 
With this setting we have
$$\E^{x_0}\int_0^\bart f(X_s)ds=\E^{x_0}\int_0^{L_\bars}f(B^1_{\tau(s)})ds.$$
where $\bars=\inf\{t:B^2(t)=0,\ |B^1(t)|\geq r\}$ and $L$ stands for the local time of $B^2$ at $0$. After the change of variable: 
$s=L_u$ in the right hand side, the equation becomes
\begin{equation}\label{eq-helper}
\E^{x_0}\int_0^\bars f(B^1_u) dL_u.
\end{equation}
Now let $g_\delta:\reals\to \reals^+$ be a family of continuous functions 
such that $\int g_\delta=1$ and $\mbox{support}(g_\delta)\subset (-\delta,\delta)$. By the occupation time formula,
we have almost surely
$$\int_0^\bars f(B^1_u) dL_u=\lim_{\de\to 0}\int_0^\bars f(B^1_u)g_\delta(B_u^2)du$$
as $\de$ tends to $0$.
Indeed a sufficient condition for this to hold is that $\bars<\infty$ almost surely. This last fact can be proved as follows: let 
$\bars_s=\inf\{t:B_t^2=0,|B_t^1|\geq s\}$. By scaling, $\forall s,s'\ \PPP(\bars_s=\infty)=\PPP(\bars_{s'}=\infty).$
So if $\bars^*=\inf\{t:B_t^2=0,\;|B^1_t|>0\}$, 
we have $\{\bars^*=\infty\}=\cap_n\{\bars_{1/n}=\infty\}.$ But $\{\bars_{1/n}=\infty\}$ is a sequence of decreasing events. So $\PPP\{\bars^*=\infty\}=
\lim_n\downarrow\PPP\{\bars_{1/n}=\infty\}=\PPP\{\bars_1=\infty\}.$ But it is easily seen that $\PPP\{\bars^*=\infty\}$ is $0$, which proves our claim.
We would then like to prove
\begin{lemma}\label{lemma-exp_cv}
\begin{equation}\label{converges?}
\E^{x_0}\int_0^\bars f(B^1_u) dL_u=\lim_{\de\to 0}\E^{x_0}\int_0^\bars f(B^1_u)g_\delta(B_u^2)du
\end{equation}
\end{lemma}
Let us postpone the proof of this Lemma, and continue with the proof of Lemma \ref{lemma-meanoccup}. 
We first rewrite the quantity in the right hand side of (\ref{converges?}). We know it is equal to $u_\delta(x_0)$ where $u_\delta$ is defined as follows:
\begin{lemma}
\label{lemma-sol_to_dir}
The unique solution to the partial differential equation:
\begin{displaymath}
\left\{\begin{array}{ll}
-1/2\Delta u_\delta(x_1,x_2)=g_\delta(x_2)f(x_1) &\textrm{on } D_{r}\\
u_\delta(x_1,x_2)=0 &\mbox{on }\partial D_{r}
\end{array}\right.
\end{displaymath}
where $D_{r}:=\reals^2 \setminus (-\infty,r]\cup[r,\infty)$
is given by 
\begin{equation}\label{eq:udelta}
u_\delta(x)=\int_{D_{r}} 2g_\delta(z_2)f(z_1)G(x,z)dz.
\end{equation}
where $G$ is the Green function of $D_{r}$. $G$ is given by a complex analogue of (\ref{eq:G}), i.e. 
\begin{equation}\label{eq:G2}
G(z_0,z)=-(1/2\pi)\log\left|\frac{h(z/r)-h(z_0/r)}{1-h(z/r)\overline{h(z_0/r)}}\right|,
\end{equation}
where, as in (\ref{eq:H}),
\begin{equation}
h(z)=\frac{\sqrt{\frac{1+z}{1-z}}-1}{\sqrt{\frac{1+z}{1-z}}+1}.
\end{equation}
\end{lemma}

\noindent
{\bf Proof of Lemma~\ref{lemma-sol_to_dir}:}
$h(z)$ can be written $u(v(z))$ where 
$$v(z):=\sqrt{\frac{1+z}{1-z}},\quad\mbox{and}$$
$$u(z):=\frac{z-1}{z+1}.$$
$v$ is a conformal mapping of $D_1$ to the upper half-plane, and $u$ is a conformal mapping of the upper-half plane to the unit disk. 
Thus $h$ maps $D_1$ conformally to the unit disk, and sends $0$ to $0$. Hence for any $z_0\in D_{r}$, 
$$\frac{h(z/r)-h(z_0/r)}{1-h(z/r)\overline{h(z_0/r)}}$$
is a conformal mapping of $D_{r}$ to the unit disk, which sends $z_0$ to $0$. Now the Green function of the unit disk with pole at $0$ is
$-(1/2\pi)\log z$. Since Green functions are conformally invariant (e.g. see \cite[p 257]{Ahlfors}), $G(z_0,z)$ as defined in (\ref{eq:G2}) is
indeed the Green function of $D_{r}$ with pole at $z_0$. Then (\ref{eq:udelta}) is simply the Green's representation formula (e.g. see
\cite[chapter 2]{GT}).\qed

So the right hand side of (\ref{converges?}) becomes 
$$u_\delta(x_0)=2\int_{-\infty}^{\infty}g_\delta(z_2)\left(\int_{-r}^{r}f(z_1)G(x_0,(z_1,z_2))dz_1\right)dz_2.$$
Now by dominated convergence, the second integral is a continuous function of $z_2$. 
Indeed, since $h$ is $1-1$ and analytic, for $z$ in a compact subset $K$ of $\Omega_{r_3}$ there exists $M$ such that 
$$\left|\frac{h(z/r)-h(x_0/r)}{1-h(z/r)\overline{h(z_0/r)}}\right| \geq M|z-x_0|.$$ Thus 
\begin{equation}\label{eq:contintmiddle}
G(x_0,z)\leq
-(1/2\pi)\log(M)-(1/2\pi)\log|x_{0,1}-z_1|
\end{equation}
which does not depend on $z_2$ and is integrable as a function of $z_1$. Here $x_{0,1}$ denotes the real part of $x_0$.
We have proved that 
$$\lim_{\delta\to 0}\E^{x_0}\int_0^\bars f(B^1_u)g_\delta(B_u^2)du = 2\int_{-r}^{r}f(z_1)G(x_0,(z_1,0))dz_1.$$
This, together with (\ref{eq-helper}) and Lemma \ref{lemma-exp_cv} yields
$$\E^{x_0}\int_0^{\bart}f(X_s)ds = 2\int_{-r_1}^{r_1}f(z_1)G(x_0,(z_1,0))dz_1.$$\qed

\noindent
{\bf Proof of Lemma~\ref{lemma-exp_cv}:}
By the occupation time formula, and the continuity of the local time,
it holds almost surely that 
\begin{equation}\label{eq-pointwisecv}
\int_0^\bars dL_u=\lim_{\delta\to 0}\int_0^\bars g_{\delta}(B_u^2)du.\end{equation}
Thus, (\ref{converges?}) would follow by applying Fatou's lemma to 
$$\int_0^\bars f(B_u^1)g_{\delta}(B_u^2)du\pm(\sup |f|)\int_0^\bars g_{\delta}(B_u^2)du$$
if we could prove 
\begin{equation}\label{dominating-eqBis}
\E^z(L_\bars)=\lim_{\delta\to 0}\E^{z}\int_0^\bars g_\delta(B^2_u)du,
\end{equation}
i.e.
\begin{equation}\label{dominating-eq}
\E^z(L_\bars)=\lim_{\delta\to 0}\int_{-\infty}^\infty g_\delta (y)\E^z(L_\bars^y)dy.
\end{equation}
where $L^y_.$ denote the local time of $B_2$ at $y$.
This last equation will hold as soon as we can show that $\E^z(L_\bars^y)$ is finite and continuous at $0$.
By (\ref{eq-pointwisecv}) and Fatou's Lemma ,
\begin{eqnarray}
\E^{z}(L^0_{\bars})&\leq&\liminf\E^z\int_0^\bars g_\delta(B^2_u)du\nn\\
&=&\liminf\int_{-\infty}^\infty\E^z(L_\bars^y)g_\delta(y)dy\label{eq-explanation1}\\
&=&\liminf \int_{-\delta}^\delta g_\delta(y)\int_{-\infty}^{\infty} G(z,(x,y))dxdy\label{eq-Fatou}
\end{eqnarray}
where (\ref{eq-explanation1}) and (\ref{eq-Fatou}) respectively follow from the occupation time formula and the Poisson representation formula of 
Lemma \ref{lemma-sol_to_dir}.
A careful study of the second integral on the right hand side of (\ref{eq-Fatou}) reveals that it is finite, and continuous as a function of $y$. 
Indeed, let us fix $\eps>0$ and define 
\begin{eqnarray*}
\lefteqn{\int_{-\infty}^{\infty}G(z,(x,y))dxdy=\int_{-\infty}^{-r_3-\eps}G(z,(x,y))dxdy}\\
& &+\int_{-r_3-\eps}^{-r_3+\eps}G(z,(x,y))dxdy+\int_{-r_3+\eps}^{r_3-\eps}G(z,(x,y))dxdy\\
& &+\int_{r_3-\eps}^{r_3+\eps}G(z,(x,y))dxdy+\int_{r_3+\eps}^{\infty}G(z,(x,y))dxdy\\
&=&I_1+I_2+I_3+I_4+I_5\\
\end{eqnarray*}
By choosing $\eps$ accordingly, $I_2$ and $I_4$ can be made arbitrarily small. And by equation (\ref{eq:contintmiddle}), $I_3$ is a continuous function 
of $y$. Finally, on $|x|\geq r_3+\eps$, $G(z,(x,y))$ can be seen to be dominated by $C/x^2$ for some constant $C>0$ uniformly on $y$ ($y$ small enough).
Thus $I_1+I_6$ tends to $0$ as $y\to 0$. These facts put together prove the continuity and the finiteness 
of the second integral in (\ref{eq-Fatou}). Therefore the 
right hand side of that equation is equal to 
$$\int_{-r_3}^{r_3}G(z,(x,0))dx$$ and in particular is finite.
Thus $\E^x(L_\bars^0)$, too, is finite. We would now like to prove that $\E^x(L_\bars^a)$ (as a function of $a$) is continuous
at $0$. Let $a>0$; then
Tanaka's formula (see e.g. \cite{Revuz-Yor}, p222) gives
\begin{eqnarray}
\frac{1}{2}L_\bars^a &=& (B^2_\bars-a)^+ -(B^2_0-a)^+ -\int_0^\bars 1_{(B^2_s>a)}dB^2_s\nn\\
&=&-\int_0^\bars 1_{(B^2_s>a)}dB^2_s.\nn
\end{eqnarray}
Therefore
$$\frac{1}{2}(L_\bars^a-L_\bars^0)=\int_0^\bars 1_{(0<B^2_s\leq a)}dB^2_s.$$
By $L^2$-isometry, we obtain
\begin{eqnarray}
\E^z\left[(\frac{1}{2}(L_\bars^a-L_\bars^0))^2\right]&=&\E^z\int_0^\bars 1_{(0<B_s^2\leq a)}ds\nn\\
&=&\int_0^a\int_{-\infty}^\infty G(z,(x,y))dxdy
\end{eqnarray}
But our previous study shows that this is finite, and tends to $0$ as $a\to 0$. So in particular that $\E^z(L_\bars^a)$ as a function 
of $a$ is continuous at $0$.
We have proved (\ref{dominating-eq}), and thus Lemma \ref{lemma-exp_cv}.\qed

\section{Upper bounds}
\label{sec-upperbound}
The following lemma will be used in proving both the lower bound and the upper bound.
 Throughout this section, fix $0<r_1\leq r_3/2$, let $X$ be a Cauchy process, $I(x,\ep)$ the open interval of radius $\ep$ centered at $x$,
$\bart:=\inf\{t>0:|X_t|\geq r_3\}$ and define 
$$\mu_{\bart}^X=\int_0^{\bart}1_{I(0,r_1)}(X_s)ds.$$

\begin{lemma}\label{lemma-meanoccup}
There exists $c>0$ such that for all $r_1\leq r_3/2$ and $|x_0|=r_2$ we have
\begin{equation}
\label{bound-meanoccup}
\E^{x_0}(\mu_{\bart}^X)\leq r_1[c+\frac{2}{\pi}\log (r_3/r_1)].
\end{equation}
and for all $k\geq 0$
\begin{equation}\label{bound-momentsmeanoccup}
\E^{x_0}(\mu_{\bart}^X)^k\leq k!r_1^k[c+\frac{2}{\pi}\log(r_3/r_1)]^k. 
\end{equation}
\end{lemma}

\noindent
{\bf Proof of Lemma~\ref{lemma-meanoccup}:}
Recall that by Proposition \ref{prop-GreenForCauchy},
$$\E^{x_0}\int_0^{\bart}1_{|X_s|\leq r_1} ds = 2\int_{-r_1}^{r_1}G(x_0,(z_1,0))dz_1.$$
To prove Lemma \ref{lemma-meanoccup} we therefore need to find an upper bound for the right hand side of the above equation, when $x_0$ lies on
the real axis. We have
\begin{eqnarray*}
\lefteqn{2\int_{-r_1}^{r_1}G(x_0,(z_1,0))dz_1 
=-\int_{-r_1}^{r_1}\frac{1}{\pi}\log \left|\frac{h(z_1/r_3)-h(x_{0,1}/r_3)}{1-h(z_1/r_3)\overline{h(x_{0,1}/r_3)}}\right|dz_1}\\
&\leq& -\int_{-r_1}^{r_1}\frac{1}{\pi}\log\left|h(z_1/r_3)-h(x_{0,1}/r_3)\right|dz_1+\frac{2}{\pi}r_1\log 2.
\end{eqnarray*}
By scaling, it suffices to work with $r_1=1$. Then by the assumption $r_3\geq 2r_1$ (see Lemma \ref{lemma-meanoccup}) we know $r_3\geq 2$. So 
we will treat the cases $|x_{0,1}|\leq (3/4)r_3$ and $|x_{0,1}|\geq (3/4)r_3$ independently. 
The function $h$, when restricted to the compact set $[-3/4,3/4]$ is smooth and its derivative does not cancel. Thus 
$$\left|h(z_1/r_3)-h(x_{0,1}/r_3)\right|\geq\left|\frac{z_1-x_{0,1}}{r_3}\right|M$$
for some constant $M>0$. Hence 
\begin{eqnarray}
-\int_{-1}^1\frac{1}{\pi}\log\left|h(z_1/r_3)-h(x_{0,1}/r_3)\right|dz_1 &\leq& -\int_{-1}^{1}{\frac{1}{\pi}}\log\left|\frac{z_1-x_{0,1}}{r_3}\right|dz_1\nn\\
& &-\int_{-1}^1\frac{1}{\pi}\log M dz_1\nn\\
&\leq& \frac{2}{\pi}\log r_3 +A\nn
\end{eqnarray}
for some $A$ uniformly in $|x_{0,1}|\leq (3/4)r_3$. It remains to treat the case $|x_{0,1}|\geq (3/4)r_3$. But then $|z_1/r_3|\leq 1/r_3\leq 1/2$
while $|x_{0,1}/r_3|\geq 3/4$. Therefore, $h$ being continuous and $1-1$, there exists a constant $B>0$ uniform in $r_3$ such that when $x_{0,1}\geq (3/4)r_3$,
$|h(z_1/r_3)-h(x_{0,1}/r_3)|\geq B$. So that in that case 
$$-\int_{-1}^1\frac{1}{\pi}\log\left|h(z_1/r_3)-h(x_{0,1}/r_3)\right|dz_1\leq C$$
for some constant $C$ uniformly in $r_3$. This finishes the proof of equation (\ref{bound-meanoccup}). (\ref{bound-momentsmeanoccup}) will then follow
from the strong Markov property for the Cauchy process. Indeed, taking $r_1=1$,
\begin{eqnarray*}
\lefteqn{\E^{x_0}(\mu_{\bart}^X)^k =k!\E^{x_0}\left(\int_{0\leq s_1\dots\leq s_k\leq\bart}\prod_{i=1}^k 1_{I(0,1)}(X_{s_i})ds_1\dots ds_k\right)}\\
&\leq& k!\E^{x_0}\left(\int_{0\leq s_1\leq \dots\leq s_{k-1}\leq\bart}\prod_{i=1}^{k-1}1_{I(0,1)}(X_{s_i})(c+{\frac{2}{\pi}}\log r_3)ds_1
\dots ds_{k-1}\right)\\
&=& k(c+\log r_3)\E^{x_0}(\mu_{\bart}^X)^{k-1},
\end{eqnarray*}
proving (\ref{bound-momentsmeanoccup}) (for $r_1=1$) by induction on $k$. Then the result for general $r_1$ follows by scaling.\qed

This leads us to
\begin{lemma}
\label{lemma-expCheb}
With the notations of Lemma \ref{lemma-meanoccup}:\\
for $0\leq \lambda<[(2/\pi)r_1\log(r_3/r_1)+cr_1]^{-1}$,
\begin{equation}\label{eq-exp}
\E(e^{\lambda\mu_\bart})\leq (1-\lambda r_1[\frac{2}{\pi}\log(r_3/r_1)+c])^{-1},
\end{equation}
implying that for $t>0$
\begin{equation}
\label{eq-Cheb}
\PPP(\mu_\bart\geq t)\leq t(r_1[\frac{2}{\pi}\log(r_3/r_1)+c])^{-1}\exp(1-t (r_1[\frac{2}{\pi}\log(r_3/r_1)+c])^{-1}).
\end{equation}
\end{lemma}

\noindent
{\bf Proof of Lemma~\ref{lemma-expCheb}:}
(\ref{eq-Cheb}) follows from (\ref{eq-exp}) by Chebychev's inequality. (\ref{eq-exp}) is a straightforward consequence of 
(\ref{bound-momentsmeanoccup}).\qed

In the remainder of this section, we use Lemma \ref{lemma-meanoccup} to prove the upper bounds in Theorem \ref{theo-sup} and Theorem \ref{theo-thick}.
Namely if we define
$$\thi_{\geq a}:=\{x\in I(0,1): \lim_{\ep\to 0}
\frac{\mu_\bart^X(I(x,\ep))}{\ep (\log \ep)^2}\geq a\}$$
(where $\bart=\inf\{t:\, |X_t|\geq 1\}$), then we will show that for any $a\in (0,2/\pi]$,\label{statement:upperthick}
\begin{equation}
\label{eq-upperthick}
\dim(\thi_{\geq a})\leq1-a\pi/2\,,\quad a.s.\,,
\end{equation}

and
\begin{equation}
\label{eq-uppersup}
\limsup_{\eps\to0}
\sup_{|x| < 1}
\frac{
\mu_\bart (I(x,\eps))}{\eps
(\log\eps)^2} \le  2/\pi
\,,\quad a.s.
\end{equation}

Our proof follows \cite{DPRZ4}. Set $h(\eps)=\eps|\log\eps|^2$ and 
$$
z (x,\ep) := \mu_{\bart} (I(x,\eps))/h(\eps).
$$
Fix $0<\delta <1$ and choose a sequence $\tilde\ep_n\downarrow 0$ as $n\rar
\ff$ in such a way that $\tilde\ep_n<e^{-2}$ and 
\begin{equation}\label{eq-geomheps}
h(\tilde\ep_{n+1}) = (1-\de) h(\tilde\ep_n),
\end{equation}
implying that $\ep_n$ is monotone decreasing in $n$. Since, for 
$\tilde\ep_{n+1}\leq \ep\leq \tilde\ep_n$
we have 
\begin{equation}
\label{eq-zepnep}
\quad \quad
z(x,\tilde\ep_n)= \frac{h(\tilde\ep_{n+1})}{h(\tilde\ep_{n})}
\frac{\mu_{\bart}(I(x,\tilde\ep_n))}{
h(\tilde\ep_{n+1})} \geq
(1-\de) z(x,\ep) \, ,
\end{equation}
it is easy to see that for any $a>0$,
\[
\mbox{\sf Thick}_{\geq a}
\subseteq D_{a}:=
\{x \in I(0,1) \,\Big|\,
\limsup_{n\rar\ff}z(x,\tilde\ep_n)\geq (1-\de)a\}.
\]
Let $\{ x_j : j=1,\ldots,\bar K_n\}$
denote a maximal collection of points in $I(0,1)$
such that $\inf_{\ell \neq j} |x_\ell-x_j| \geq \delta \tilde\ep_n$.
Let $\bart_2:=\inf\{t:\,|X_t|\geq 2\}$ and 
$\AA_{n}$ be the set of $1\leq j\leq \bar K_n$, such that
\begin{equation}
\label{ed-defofAn}
\mu_\bart^X(I(x_j,(1+\delta)\tilde\ep_n))\geq(1-2\delta)ah(\tilde\ep_n).
\end{equation}
Applying (\ref{eq-Cheb}) with $r_1=(1+\delta)\tilde\ep_n$ and $r_3=2$ gives
$$
\PPP^{x}(\mu_{\bart_2}(I(0,(1+\delta)))\geq(1-2\delta)a h(\tilde\ep_n))\leq c{\tilde\ep_n}^{\,\,a(1-5\delta)\pi/2},
$$
for some $c=c(\delta)<\ff,$ and any $x\in I(0,1)$. Note that for all $x\in I(0,1)$ and $\ep,b\geq 0$
$$
\PPP(\mu_\bart(I(x,\eps))\geq b)\leq \PPP^{-x}(\mu_{\bart_2}(I(0,\ep))\geq b).
$$
Thus for any $j$ and $a>0$,
$$\PPP(j\in \AA_n)\leq c\tilde\ep_n^{\,\,a(1-5\delta)\pi/2},$$
implying that
\begin{equation}\label{eq-boundexpect}
\E|\AA_n|\leq c'{\tilde\ep_n}^{\,\,a(1-5\delta)\pi/2-1}
\end{equation}
(by definition of $\bar K_n$).
Let $\VV_{n,j}=I(x_j,\de\tilde\ep_n)$. For any $x \in
D(0,1)$ there exists $j \in \{1,\ldots,\bar K_n \}$ such that
$x \in \VV_{n,j}$,
hence
$I(x,\tilde\ep_n) \subseteq I(x_j,(1+\de)\tilde\ep_n)$. Consequently,
$\cup_{n \geq m} \cup_{j\in \AA_{n}}\VV_{n,j}$ forms a
cover of
$D_{a}$
by sets of maximal diameter $2\de \tilde\ep_m$.
Fix $a \in (0,2/\pi]$.
Since $\VV_{n,j}$ have diameter $2\de\tilde\ep_n$,
it follows from (\ref{eq-boundexpect}) that for $\gamma=1-\pi a(1-6\de)/2  > 0$,
\[
\E \sum_{n=m}^\infty \sum_{j\in \AA_{n} }|\VV_{n,j}|^\gamma \leq
c'\, (2\de)^\gamma \sum_{n=m}^\infty \tilde\ep_{n}^{\,\,\de a\pi/2 }
 < \infty \;.
\]
Thus, $\sum_{n=m}^\infty \sum_{j\in \AA_n}
|\VV_{n,j}|^\gamma$ is
 finite
a.s. implying that $\dim (D_{a})\leq \gamma$  a.s.
Taking $\de \downarrow 0$
completes the proof of the upper bound (\ref{eq-upperthick}).

Turning to prove (\ref{eq-uppersup}), set $a=2(1+\de)/(\pi(1-5\de))$ noting that
by (\ref{eq-boundexpect})
$$ \sum_{n=1}^\infty \PPP(|\AA_n|\geq 1)
\leq \sum_{n=1}^\infty \E|\AA_n| \leq c'\sum_{n=1}^\infty
\tilde\ep_n^{\,\,\de}< \infty\,.$$
By Borel-Cantelli, it follows that
a.s. $\AA_n$ is empty for all
$n>n_0(\omega)$ and some $n_0(\om)<\infty$. By (\ref{eq-zepnep}) we then have
$$
\sup_{\ep \leq \tilde\ep_{n_0(\om)}} \, \sup_{|x|<1}
\frac{\mu_{\bart} (I(x,\eps))}{\eps
(\log\eps)^2} \le  a\frac{1-2\de}{1-\de}\leq a\,,\,$$
and (\ref{eq-uppersup}) follows by taking $\delta\downarrow 0$.
\qed

\section{Lower bounds}

In this section we adapt the proof of \cite[section 3]{DPRZ5}. While its authors studied the intersection local time for two independent Brownian motions,
we are interested in the same quantity but for the intersection of a Brownian motion with a line. Throughout what follows we use notations similar to
those in \cite[section 3]{DPRZ5}.

Fixing $a<2/\pi$, $c>0$ and $\delta>0$,
let 
\[
\bart_c:=\inf\{t>0:\;|X_t|\geq c\}\;,
\]
\[
\Gamma_c = \Gamma_c(\om) := \{x\in I(0,1):\;
\lim_{\eps\to 0} \frac{\mu_{\bart_c} (I(x,\eps))}{\eps
(\log\eps)^2}= a\} \;,
\]

 and  ${\mathcal E}_c := \{ \om : \dim(\Gamma_c(\om) ) \geq 1-a\pi/2
-\delta \}.$

In view of the results of Section \ref{sec-upperbound}, we will obtain Theorem
 \ref{theo-thick} once we show that
$\PPP({\mathcal E}_1)=1$ for any $a<1$ and $\de>0$. Indeed, since $\thi_a\subset\thi_{\geq a}$ and since we have seen that 
$$\dim(\thi_{\geq a})\leq 1-a\pi/2 \,,\quad a.s.,$$ proving $\PPP({\mathcal E}_1)=1$ will imply $$\dim \thi_a = 1-a\pi/2 \,,\quad a.s.$$
Moreover, then the inequality
$$
\liminf_{\ep\to 0} \sup_{|x|<1}
\frac{\mu_{\bart} (I(x,\eps))}{\eps
(\log \eps)^2}
\geq
\sup_{|x|<1} \liminf_{\ep\to 0}
\frac{\mu_{\bart} (I(x,\eps))}{\eps
(\log\eps)^2}
$$
implies that for any $\eta>0$,
$$
\liminf_{\eps\to0}
\sup_{|x|<1}
\frac{\mu_{\bart} (I(x,\eps))}{\eps
(\log\eps)^2} \ge  2/\pi-\eta
\,,\quad a.s.
$$
In view of (\ref{eq-uppersup}), these lower bounds establish Theorem \ref{theo-sup}.

The bulk of this section and the next will be dedicated to showing that
$\PPP({\mathcal E}_1)>0$. Assuming
this for the moment, let us show that this implies
$\PPP({\mathcal E}_1)=1$.
With $X^c_t:=c^{-1}X_{ct}$ we have that
$c\bart(\om^c)= \inf\{ct:\,|c^{-1}X_{ct}| \geq 1\}=\bart_c(\om)$, and hence
\begin{eqnarray}
\label{eq-indepofc}
\mu^{X^c}_{\bart} (I(x,\eps)) &&= \,\,\int_0^{\bart(\om^c)}1_{\{|X^c_s-x|\leq\ep\}}\,ds\,
=\,\int_0^{\bart(\om^c)}1_{\{|X_{cs}-cx|\leq 
c\ep\}}\,ds\nn\\
&&=\frac{1}{c}\int_0^{c\bart(\om^c)}1_{\{|X_{s}-cx|\leq
c\ep\}}\,ds\nn \,=\, \frac{1}{c}
\mu^X_{\bart_c} (I(cx,c\eps)). \nn
\end{eqnarray}
Consequently, $\Gamma_c(\om)=c\Gamma_1(\om^c)$, so the Cauchy process' scaling property implies that  $p=\PPP({\mathcal E}_c)$ is
independent of
$c>0$. Let
$$
{\mathcal E} := \limsup_{n \to \infty} {\mathcal E}_{n^{-1}} \;,
$$
so that $\PPP({\mathcal E}) \geq p$. The Cauchy process is a Feller process; hence if we let $({\mathcal F}_t)$ be the usual augmentation 
of the natural filtration, it can be shown that $({\mathcal F}_t)$ is right-continuous (e.g. see \cite[III-2]{Revuz-Yor}). Therefore, since
${\mathcal E}_c \in {\mathcal F}_{\bart_c}$, ${\mathcal E}\in {\mathcal F}_0$ which implies $\PPP({\mathcal E})\in\{0,1\}$. Thus, $p>0$ yields
$\PPP(\mathcal{E})=1$. We will see momentarily that the events $\mathcal{E}_c$ are essentially increasing in $c$, i.e.
\begin{equation}\label{eq-essinc}
\forall\, 0<b<c\quad\PPP(\mathcal{E}_b\setminus \mathcal{E}_c)=0.
\end{equation}
Thus, $\PPP({\mathcal E}\setminus {\mathcal E}_1) \leq
\PPP(\bigcup_{n} \{
{\mathcal E}_{n^{-1}} \setminus {\mathcal E}_1 \}
) = 0$, so
that also $\PPP({\mathcal E}_1)=1$.
To see (\ref{eq-essinc}), we proceed exactly as in \cite{DPRZ4}. First notice that for $b<c$,
$$\Gamma_b(\om)\setminus\{\om:\,0\leq t\leq\bart_c\}\subset\Gamma_c(\om).$$
Hence
$$\PPP(\mathcal{E}\setminus\mathcal{E}_c)\leq \E\PPP(\dim(\Gamma_b(\om))\neq \dim (\Gamma_b(\om)\setminus\{\om:\,\bart_b\leq t\leq\bart_c\})|
\mathcal{F}_{\bart_b}).$$
Then applying the strong Markov property at time $\bart_b$ and observing that the set $\Gamma_b(\om)$ is Borel gives 
(\ref{eq-essinc}) exactly as in \cite{DPRZ4} (since the Cauchy process does not hit points).

So we just have to show that $\PPP(\mathcal{E}_1)>0$. To achieve this goal we will use the Brownian representation of the Cauchy process, and follow the
strategy of \cite{DPRZ5}. More precisely, moving to a Brownian setting, we will now focus our attention on the ``projected intersection 
local time measures'':
$$\II_{t} (A):=\int_0^t 1_{\{B_u\in A\}}dL_u^0,$$
where $B=(B^1,B^2)$ is a planar Brownian motion and $L_.^0$ is the local time of $B^2$ at $0$. 
$\II$ is simply the amount of local time spent in $A$ before $t$. To see how this relates to the Cauchy process, note that, for example, for any 
set $A\subset \R^2$, $\II_{\bars}(A)=\mu_\bart(A\cap \mbox{$x$-axis})$ where
$\bars:=\inf\{t:B^2_t=0\mbox{ and }|B^1_t|\geq 1\}$, $\bart:=\inf\{t:|X_t|\geq 1\}$, $X$ is the Cauchy process associated to the planar Brownian motion 
$B$ and $\mu$ is the occupation measure for $X$. We then reproduce the setting of \cite[p 248]{DPRZ5}:
fix $a<2$,
$\ep_1=1/8$
and the square
$S=S_1=[\epon,2\epon ]^2 \subset D(0,1)$. Note that
for all $x \in S$ and $y \in S \cup \{0\}$ both $0 \notin D(x,\epon)$ and
$0 \in D(x,1/2) \subset D(y,1) \subset D(x,2)$.
Let $\ep_k=\epon (k!)^{-3}=\epon \prod_{l=2}^k l^{-3}$.
For $x \in S$, $k \geq 2$
and $\rho>\epon$, let $N_k^x(\rho)$ denote the number of
excursions of $B_\cdot$ from $\partial D(x,\ep_{k-1})$ to
$\partial D(x,\ep_{k})$ prior to hitting $\partial D(x,\rho)$.
Set $n_k=3ak^2\log k$. We will say that a point $x\in S$ is
{\bf n-perfect} if
$$
n_k-k\leq N^x_k(1/2) \leq N^x_k(2) \leq n_k+k\,,
\hspace{.3in}\forall k=2,\ldots,n.
$$

For $n \geq 2$ we partition $S$ into
$M_n=\epon^2/(2 \ep_n)^2=(1/4) \prod_{l=1}^n l^6$ non-overlapping squares
of edge length $2 \ep_n=2 \epon/(n!)^3$, which we denote by
$S(n,i)\,;\,i=1,\ldots,M_n$ with $x_{n,i}$ denoting the center of
each $S(n,i)$. Let
$Y(n,i)\,;\,i=1,\ldots,M_n$ be the sequence of random variables defined by
\[Y(n,i)=1\hspace{.1in}\mbox{if $x_{n,i}$ is n-perfect}\]
and $Y(n,i)=0$ otherwise. 
Define
$$
A_n=\bigcup_{i: Y(n,i)=1}S(n,i),
$$
and
\beq{F-def}
F=F(\om)=\bigcap_m \overline{\bigcup_{n\ge m} A_n}
:=\bigcap_m F_m
\,.
\end{equation}
Note that each $x \in F$ is the limit of a sequence $\{x_n\}$ such that
$x_n$ is $n$-perfect. We finally rotate this picture by $45$ degrees to the right. $S$ now intersects the $x$-axis; let $D$ be this intersection.
The next lemma will be proved in the next section.
\begin{lemma}\label{lemma-thick=perfect}
Let $\bart:=\inf\{t:\;|B_t|\geq 1\}$. A.~s.~ for all $x\in F\cap D$
$$
\lim_{\eps\to 0} \frac{\II_{\bart} (D(x,\eps))}{\eps
(\log\eps)^2}= \frac{2}{\pi}a.
$$
\end{lemma}
Now Lemma $3.2$ in \cite{DPRZ5} shows that for every $a<1$, and for every $\delta >0$ such that $1-a-\delta>0$,
\begin{equation}\label{eq-rightdim}
\PPP(\dim(F\cap D)\geq 1-a-\de)>0.
\end{equation}
This, together with Lemma \ref{lemma-thick=perfect} implies
$$\PPP\left(\dim\left(D\cap \left\{\lim_{\ep\rar 0}\frac{\II_\bart (D(x,\ep))}{\ep (\log\ep)^2}=\frac{2}{\pi}a\right\}\right)\geq
1-a-\de\right) >0.$$
Now if $\bars$ denotes the first time that the planar Brownian motion hits the complement of $(-1,1)$
 (on the real axis),
 we have by the strong Markov property
\begin{multline*}
\PPP\left(\dim\left(D\cap \left\{\lim_{\ep\rar 0}\frac{\II_\bars (D(x,\ep))}{\ep (\log\ep)^2}=\frac{2}{\pi}a\right\}\right)\geq
1-a-\de\right)\\
\geq K\PPP\left(\dim\left(D\cap \left\{\lim_{\ep\rar 0}\frac{\II_\bart (D(x,\ep))}{\ep (\log\ep)^2}=\frac{2}{\pi}a\right\}\right)\geq
1-a-\de\right)
\end{multline*}
where 
$$K=\inf_{|a|=1}\PPP^a\left(\mbox{Brownian motion hits }(-\ff,-1]\cup[1,\ff)\mbox{ before D}\right)$$
Since $K>0$, we have proved that $\PPP(\mathcal{E}_1)>0$, which concludes the section.
\section{From excursions to intersection local time}
This section follows closely the argument developed in \cite[section 4]{DPRZ5}.
The sets $F$ and $D$ are the same as in the previous section, and 
$h(\ep):=\ep(\log\ep)^2$. Lemma \ref{lemma-thick=perfect} will follow from the next two lemmas.
\begin{lemma}\label{lemma-upthick}
For every $\de>0$, if $x\in F\cap D$ then 
\begin{equation}
\frac{2}{\pi}a(1-\de)^5\leq\liminf_{\ep\rar 0}\II_{\bart}(D(x,\ep))/h(\ep)
\end{equation}
\end{lemma}
\begin{lemma}\label{lemma-lowthick}
For every $\de>0$, if $x\in F\cap D$ then 
\begin{equation}
\limsup_{\ep\rar 0}\II_{\bart}(D(x,\ep))/h(\ep)\leq \frac{2}{\pi}a(1+\delta)^5.
\end{equation}
\end{lemma}

\noindent
{\bf Proof of Lemma~\ref{lemma-upthick}:} We use the same notations as in \cite{DPRZ5}. Let $\de_k=\ep_k/k^6$ and let $\mathcal{D}_k$ be a 
$\de_k$-net of points in
$S$. Let
\[\ep'_k=\ep_ke^{1/k^{6}},\hspace{.2in}\ep''_{k-1}=\ep_{k-1}e^{-1/k^{6}},\]
so that
\begin{equation}
\ep'_k\geq\ep_k+\de_k,\hspace{.2in}\ep''_{k-1}\leq \ep_{k-1}-\de_k.\label{21.1}
\end{equation}
We will say that a point $x'\in \mathcal{D}_k$ is
{\bf lower k-successful} if there
are at least
$n_k-k$ excursions of $W$ from $\partial D(x',\ep''_{k-1})$ to
$\partial D(x',\ep'_k)$
 prior to
$\bar{\th}$.
Let \[\ep_{k,j}=\ep_ke^{-j/k}\,,\,j=0,1,\ldots,3k\log (k+1),\]
and let $\ep'_{k,j}=\ep_{k,j}e^{-4/k^{3}}=\ep'_ke^{-j/k}e^{-4/k^{3}-1/k^{6}}$
(the choice of $\ep'_{k,j} =\ep_{k,j}e^{-2/k^{3}}$ used in \cite[Section 6]{DPRZ4} and \cite[Section 5]{DPRZ5}
works well in the former, but is neither
suitable in the context of the latter, nor here, since then
$\ep'_{k,j} + 2\ep_k/k^6\geq \ep_{k,j} >\ep_{k,j}$ for large $k$ whereas the reverse inequality is needed in \cite[page 4984]{DPRZ5}.
Taking $\ep_{k,j}= \ep_{k,j}e^{-4/k^{3}}$ as done here fixes this problem.)
By analogy with \cite{DPRZ5}, we say that $x'\in \mathcal{D}_k$ is {\bf lower k,$\delta$-successful} if it is
lower k-successful and in addition,
\begin{equation}
\frac{2}{\pi}(1-\de) \ep'_{k,j}\le \rho
(D(x',\ep'_{k,j})) ,
\;\forall j=0,\ldots,3k\log (k+1).
\end{equation}
where $\rho$ denotes the measure supported on the real axis, and whose restriction to the real axis is $1/\pi$ times the Lebesgue measure.
We recall Lemma 2.3 of \cite{DPRZ5}, adapted to our situation. In what follows , $\bart_{x,r}:=\inf\{t:\, |B_t-x|=r\}$.

\begin{lemma}
\label{lemma-exprho}

We can find $c<\ff$ such that for all $k \geq 1$,
$r_1\leq r_2\leq
r/2\leq
1/2
$,
$x$ and
$x_0$ with
$|x_0-x|=r_2$,
\begin{equation}
\label{eq-secmo}
\E^{x_0}(\II_{\bart_{x,r}}(D(x,r_1)))^k\leq k!
\( \rho(D(x,r_1)) \log (r/r_1)
+ c r_1  \)^k \,,
\end{equation}
and
\begin{equation}
\label{eq-+-}
\E^{x_0}(\II_{\bart_{x,r}}(D(x,r_1)))=
\rho(D(x,r_1)) \log (r/r_2) \pm c r_1 \,.
\end{equation}
\end{lemma}

The above lemma can be seen as an analogue of Lemma \ref{lemma-meanoccup}; the main difference lies in the fact that the Brownian motion is now stopped
when it leaves a disk of radius $r$. We are now in a position to prove Lemmas \ref{lemma-upthick} and \ref{lemma-lowthick}. We will derive
Lemma \ref{lemma-upthick} from the following lemma, which is an analogue of Lemma 4.3 in \cite{DPRZ5}.
\begin{lemma}
\label{lemma-lowersuc}
There exists a $k_0=k_0(\de,\om)$ such that for all $k\geq k_0$ and $x'\in\mathcal{D}_k$, if $x'$ is lower $k,\de$-successful then
\begin{equation}
\frac{2}{\pi}(1-\de)^4 h(\ep'_{k,j})\leq \II_\bart(D(x',\ep'_{k,j})), \hspace{.2in}\forall j=0,1,\ldots,3k\log (k+1).
\end{equation}
\end{lemma}

The derivation of Lemma \ref{lemma-upthick} from Lemma \ref{lemma-lowersuc} 
is exactly the same as in \cite{DPRZ5}, except that one should use $\rho$ instead of $\mu_\bart^{W'}$ when
verifying that for $k$ large enough, $x_k$ is lower $k,\de$-successful.

\noindent{\bf Proof of Lemma \ref{lemma-lowersuc}:}
The main difference with the proof in \cite[page 4984]{DPRZ5} is that now
$$
A(x',k,j)
= \{\II_{\bart} (D(x',\ep'_{k,j}))
\le \frac{2}{\pi}a
(1-\de)^4 h(\ep'_{k,j})\}.
$$
and 
$$\rho(D(0,r_1))\geq \frac{2}{\pi}(1-\de)r_1.$$
So using Lemma \ref{lemma-exprho} one obtains
$$\E_{x',s}(\tau_{l,k,j})
\geq
\frac{2}{\pi} (1-\de)^2 \log (r/r_2) r_1 \,.
$$
This, combined with the Stirling's approximation \cite[(4.6)]{DPRZ5} yields 
$$
\frac{2}{\pi}a(1-\de)^3 h(\ep'_{k,j}) \leq n'_k \E_{x',s}(\tau_{l,k,j}) \,.
$$
The remainder of the proof is identical to the one in \cite{DPRZ5}.

\noindent{\bf Proof of Lemma \ref{lemma-lowthick}:} Again, we use the same notations as in \cite{DPRZ5}: we let
\[\bar{\ep}'_k=\ep_ke^{-2/k^{6}},\hspace{.2in}\bar{\ep}''_{k-1}=\ep_{k-1}e^{
1/k^{6}},\]
so that
$$
\bar{\ep}'_k\leq\ep_k-\de_k,\hspace{.2in}\bar{\ep}''_{k-1}\geq
\ep_{k-1}+\de_k.
$$
We now say that $x'\in \mathcal{D}_k$ is {\bf upper k-successful} if
there are at most $n_k+k$ excursions of $W$ from
$\partial D(x',\bar{\ep}''_{k-1})$ to $\partial D(x',\bar{\ep}'_k)$
prior to $\bar{\th}$.
In our case, there is no need to define upper $k,\de$-successful points since our measure $\rho$ is deterministic and $\rho(A)$ is bounded by the diameter
of $A$. Using the same arguments as in the previous case, Lemma \ref{lemma-lowthick} can be derived from
\begin{lemma}\label{lemma-uppersuc}
There exists a $k_0=k_0(\de,\om)$
such that for all $k\geq k_0$ and
$x'\in  \mathcal{D}_k$, if $x'$ is upper k-successful
then
$$
\frac{2}{\pi}a
(1+\de)^4 h(\ep'_{k,j}) \ge \II_{\bart}
(D(x',\ep'_{k,j})) ,
\hspace{.2in}\forall j=0,1,\ldots,3k\log (k+1).
$$
\end{lemma}

\noindent{\bf Proof of Lemma \ref{lemma-uppersuc}:}
Again, the proof is very similar to the one of \cite[Lemma 4.4]{DPRZ5}. The three inequalities at the bottom of \cite[page 19]{DPRZ5} become
$$
 \E_{x',us}(\tau_{l,k,j}) \leq
\frac{2}{\pi}(1+\de) \log(\bar{\ep}''_{k-1}/\bar{\ep}'_k)
 (\ep'_{k,j}) \;,
$$
$$
n_k'' \log(\bar{\ep}''_{k-1}/\bar{\ep}'_k)
\leq a (1+\de) |\log (\ep'_{k,j})|^2 \;.
$$
and
$$
\frac{2}{\pi}a (1+\de)^3 h(\ep'_{k,j}) \geq n_k'' \E_{x',us}(\tau_{l,k,j}).
$$
Finally, in our case $\wh{\tau}_{l,k,j}$ becomes
$$
\wh{\tau}_{l,k,j} := \frac{\tau_{l,k,j}}
{\log (\bar{\ep}''_{k-1}/\ep'_{k,j}) \ep^{'}_{k,j}}.
$$
These are the only differences between the two proofs.

\vspace{2mm}

\noindent{\bf Acknowledgements} I am very grateful to Amir Dembo and Yuval Peres for suggesting this problem to me, and for many helpful discussions.

\end{document}